\documentclass[12pt]{amsart}
\usepackage{amssymb}
\usepackage{amsfonts}
\usepackage{amsbsy}
\usepackage{latexsym}
\newtheorem{theorem}{Theorem}[section]
\newtheorem{lemma}[theorem]{Lemma}

\newtheorem{remark}[theorem]{Remark}
\theoremstyle{definition}

\begin{document}
\title {On Jordan Derivations of  Triangular algebras}
\author{Xuehan Cheng}
\address{College of  Mathematics and Information, Ludong University,  Yantai, 264025 P. R. China}

\author{Wu Jing}
\address{Department of Mathematics and Computer Science, Fayetteville State University,  Fayetteville, NC 28301}
 \email{wjing@uncfsu.edu}

\thanks{This work is partially supported by NNSF of China (No. 10671086) and NSF of Ludong University (No. LY20062704).}
\date{June 12, 2007}
\begin{abstract} In this short note we prove that every Jordan derivation of triangular algebras is a derivation.
  \end{abstract} \maketitle

\section{Introduction}

Suppose that $\mathcal {A}$ is an algebra over a commutative ring
$\mathcal {R}$. An $\mathcal {R}$-linear map $\delta $ from
$\mathcal {A}$ to an $\mathcal {A}$-module $\mathcal {M}$ is said to
be a \textit{derivation} if $\delta (ab)=\delta (a)b+a\delta (b)$
for any $a, b\in \mathcal {A}$. We call $\delta $ a \textit{Jordan
derivation} if $\delta (a^2)=\delta (a)a+a\delta (a)$ for each $a\in
\mathcal {A}$. Obviously, every derivation is a derivation. But the
inverse, in general, is not true (see \cite{ben}).

It is natural and very interesting to find some conditions under
which each Jordan derivation is a derivation. The first result on
this direction was due to Herstein in 1957. He proved that every
Jordan derivation on a $2$-torsion free prime ring is a derivation.
In 1988, Bre$\check{\textrm {s}}$ar generalized Herstein's result to
Jordan derivations of semiprime rings (\cite{bresar}).   Recently,
Zhang and Yu (\cite{zhang}) showed that every Jordan derivation of a
triangular algebra is a derivation. More precisely, they proved the
following result.

\begin{theorem}(\cite{zhang})\label{zhangtheorem}
Let $\mathcal {A}, \mathcal {B}$ be unital algebras over a
$2$-torsion free commutative ring $\mathcal {R}$, and $\mathcal {M}$
be a unital $(\mathcal {A}, \mathcal {B})$-bimodule that is faithful
as a left $\mathcal {A}$-module and also as a right $\mathcal
{B}$-module. Then  every Jordan derivation from the triangular
algebra $Tri(\mathcal {A}, \mathcal {M}, \mathcal {B})$ into itself
is a derivation.
\end{theorem}

Note that this result requires that both $\mathcal {A}$ and
$\mathcal {B}$ are unital, and the proof is heavily dependent on
\begin{displaymath} \mathcal P=\left( \begin{array}{ll}
1 & 0\\
0 & 0
\end{array}\right) \textrm{and} \,
 Q=\left( \begin{array}{ll}
0 & 0\\
0 & 1
\end{array}\right).
 \end{displaymath}

The aim of this short note is to generalize Theorem
\ref{zhangtheorem} to more general case. We want to mention here
that we do not require the existence of identities for both
$\mathcal {A}$ and $\mathcal {B}$ and  our approach is quite
different from that in  \cite{zhang}.

We now introduce some definitions and notations. Recall that a
\textit{triangular algebra} $\mathcal {T}=Tri(\mathcal {A}, \mathcal
{M}, \mathcal {B})$ is  an algebra of the form

 \begin{displaymath}Tri(\mathcal {A}, \mathcal {M}, \mathcal {B})=\{
\left( \begin{array}{ll}
a & m\\
0 & b
\end{array}\right): a\in \mathcal {A}, m\in \mathcal {M}, b\in \mathcal {B}\}
\end{displaymath}
under the usual matrix operations, where $\mathcal {A}$ and
$\mathcal {B}$ are two algebras over a commutative ring $\mathcal
{R}$, and $\mathcal {M}$ is an $(\mathcal {A}, \mathcal
{B})$-bimodule which is faithful as a left $\mathcal {A}$-module and
also as a right $\mathcal {B}$-module (see \cite{cheung}).

Throughout this paper, we set \begin{displaymath} \mathcal
{T}_{11}=\{ \left( \begin{array}{ll}
a & 0\\
0 & 0
\end{array}\right): a\in \mathcal {A} \},
\end{displaymath}
\begin{displaymath}\mathcal {T}_{12}=\{
\left( \begin{array}{ll}
0 & m\\
0 & 0
\end{array}\right): m\in \mathcal {M} \},
\end{displaymath}
and
\begin{displaymath} \mathcal {T}_{22}=\{
\left( \begin{array}{ll}
0 & 0\\
0 & b
\end{array}\right): b\in \mathcal {B} \}.
\end{displaymath}
Then we may write $\mathcal {T}=\mathcal {T}_{11}\oplus \mathcal
{T}_{12}\oplus \mathcal {T}_{22}$, and every element $a\in\mathcal
{T}$ can be written as $a=a_{11}+a_{12}+a_{22}$. Note that notation
$a_{ij}$ denotes an arbitrary element of $\mathcal {T}_{ij}$.

\section {Jordan derivation of triangular algebras}
Throughout this section, $\mathcal {A}$ and $\mathcal {B}$ will be
two algebras over a $2$-torsion free commutative ring $\mathcal {R}$
with the property:

 (P) Suppose that  $a\in \mathcal {A}$ (resp. $\mathcal {B}$). If $xay+yax=0$ holds for all $x, y\in \mathcal {A}$ (resp. $\mathcal {B}$), then $a=0$.

Map $\delta $ will be a Jordan derivation from triangular algebra
$\mathcal {T}=Tri(\mathcal {A}, \mathcal {M}, \mathcal {B})$ into
itself, where $\mathcal {M}$ is a faithful $(\mathcal {A}, \mathcal
{B})$-bimodule.

We begin with the following useful lemma.

\begin{lemma}\label{jj}
Let $a$ be in $\mathcal {A}$ (resp. $\mathcal {B}$).

(i) If $xax=0$ for any $x\in \mathcal {A}$ (resp. $\mathcal {B}$),
then $a=0$;

(ii) If $ax+xa=0$ for any $x\in \mathcal {A}$ (resp. $\mathcal
{B}$), then $a=0$.
\end{lemma}
\begin{proof}
(i) It follows directly from Property (P).

(ii) For arbitrary $x, y\in \mathcal {A}$, multiplying $ax+xa=0$
from the left and the right by $y$ respectively and adding them
together, we obtain
$$yax+xay+yxa+axy=0.$$
This leads to $$yax+xay-ayx-xya=0.$$ Furthermore, we have
$$yax+xay+yax+xay=0.$$
By Property (P), we see that $a=0$.
\end{proof}

The following result is given by Herstein.

\begin{lemma} Let $\delta $ be a Jordan derivation on a $2$-torsion free ring $\mathcal {R}$ into itself. For any $a, b, c\in \mathcal {R}$, the following hold.

(i) $\delta (ab+ba)=\delta (a)b+a\delta (b)+\delta (b)a+b\delta
(a)$;

(ii) $\delta (aba)=\delta (a)ba+a\delta (b)a+ab\delta (a)$;

(iii) $\delta (abc+cba)=\delta (a)bc+a\delta (b)c+ab\delta
(c)+\delta (c)ba+c\delta (b)a+cb\delta (a)$.
\end{lemma}

\begin{lemma}\label{1122}
For arbitrary $a_{11}\in \mathcal {T}_{11}$ and $b_{22}\in \mathcal
{T}_{22}$, we have

(i) $\delta (a_{11})_{22}=0$;

(ii) $\delta (b_{22})_{11}=0$;

(iii) $\delta (a_{11}b_{22})=\delta (a_{11}b_{22})+a_{11}\delta
(b_{22})$;

(iv) $\delta (b_{22}a_{11})=\delta (b_{22})a_{11}+b_{22}\delta
(a_{11})$.
\end{lemma}
\begin{proof}
We compute
\begin{eqnarray*}
0&=&\delta (a_{11}b_{22}+b_{22}a_{11})\\
&=&\delta (a_{11})b_{22}+a_{11}\delta (b_{22})+\delta (b_{22})a_{11}+b_{22}\delta (a_{11})\\
&=&\delta (a_{11})_{12}b_{22}+\delta (a_{11})_{22}b_{22}+a_{11}\delta (b_{22})_{11}\\
& &+a_{11}\delta (b_{22})_{12}+\delta
(b_{22})_{11}a_{11}+b_{22}\delta (a_{11})_{22}.
\end{eqnarray*}
It follows that $$\delta (a_{11})_{22}b_{22}+b_{22}\delta
(a_{11})_{22}=0$$ and
$$a_{11}\delta (b_{22})_{11}+\delta (b_{22})_{11}a_{11}=0.$$

By Lemma \ref{jj}, we have $\delta (a_{11})_{22}=\delta
(b_{22})_{11}=0$.

Now (iii) and (iv) follows easily.
\end{proof}

\begin{lemma} \label{12} For any $a_{12}\in \mathcal {T}_{12}$, the following are true.

(i) $\delta (a_{12})_{11}=0$;

(ii) $\delta (a_{12})_{22}=0$.
\end{lemma}
\begin{proof}
For arbitrary $b_{11}\in \mathcal {T}_{11}$, we have
\begin{eqnarray*}
0&=&\delta (b_{11}a_{12}b_{11})\\
&=&\delta (b_{11})a_{12}b_{11}+b_{11}\delta (a_{12})b_{11}+b_{11}a_{12}\delta (b_{11})\\
&=&b_{11}\delta (a_{12})_{11}b_{11}+b_{11}a_{12}\delta (b_{11})_{22}\\
&=&b_{11}\delta (a_{12})_{11}b_{11}.
\end{eqnarray*}
It follows from Lemma \ref{jj} that $\delta (a_{12})_{11}=0$.

Similarly, we can get $\delta (a_{12})_{22}=0$ by considering
$\delta (b_{22}a_{12}b_{22})$ for any $b_{22}\in \mathcal {T}_{22}$.
\end{proof}
\begin{lemma} For any $a_{11}\in \mathcal {T}_{11}$ and $b_{12}\in \mathcal {T}_{12}$,

(i) $\delta (b_{12}a_{11})=\delta (b_{12})a_{11}+b_{12}\delta
(a_{11})$;

(ii) $\delta (a_{11}b_{12})=\delta (a_{11})b_{12}+a_{11}\delta
(b_{12})$.
\end{lemma}
\begin{proof}
(i) Note that
\begin{eqnarray*}
\delta (b_{12})a_{11}+b_{12}\delta (a_{11})&=&\delta (b_{12})_{11}a_{11}+b_{12}\delta (a_{11})_{22}\\
&=&0=\delta (b_{12}a_{11}).
\end{eqnarray*}

(ii) We have
\begin{eqnarray*}
\delta (a_{11}b_{12})&=&\delta (a_{11}b_{12}+b_{12}a_{11})\\
&=&\delta (a_{11})b_{12}+a_{11}\delta (b_{12})+\delta (b_{12})a_{11}+b_{12}\delta (a_{11})\\
&=&\delta (a_{11})b_{12}+a_{11}\delta (b_{12}).
\end{eqnarray*}
\end{proof}

Similarly, we have
\begin{lemma} For arbitrary $a_{12}\in \mathcal {T}_{12}$ and $b_{22}\in \mathcal {T}_{22}$,

(i) $\delta (a_{12}b_{22})=\delta (a_{12})b_{22}+a_{12}\delta
(b_{22})$;

(ii) $\delta (b_{22}a_{12})=\delta (b_{22})a_{12}+b_{12}\delta
(a_{12})$.
\end{lemma}

\begin{lemma} $\delta $ is a derivation on $\mathcal {T}_{12}$.
\end{lemma}
\begin{proof}
For any $a_{12}, b_{12}\in \mathcal {T}_{12}$, by Lemma \ref{12}, we
have
\begin{eqnarray*}
\delta (a_{12})b_{12}+a_{12}\delta (b_{12})&=&\delta (a_{12})_{11}b_{12}+a_{12}\delta (b_{12})_{22}\\
&=&0=\delta (a_{12}b_{12}).
\end{eqnarray*}
\end{proof}

\begin{lemma}
$\delta $ is a derivation on $\mathcal {T}_{11}$.
\end{lemma}
\begin{proof}
Let $a_{11}, b_{11}\in \mathcal {T}_{11}$, and $c_{12}\in \mathcal
{T}_{12}$ be arbitrary. On one side, we have
\begin{eqnarray*}
& &\delta (a_{11}b_{11}c_{12})\\
&=&\delta (a_{11})b_{11}c_{12}+a_{11}\delta (b_{11}c_{12})\\
&=&\delta (a_{11})b_{11}c_{12}+a_{11}\delta
(b_{11})c_{12}+a_{11}b_{11}\delta (c_{12}).
\end{eqnarray*}
On the other side, we get $$\delta (a_{11}b_{11}c_{12})=\delta
(a_{11}b_{11})c_{12}+a_{11}b_{11}\delta (c_{12}).$$ Then we can
infer that
$$\delta (a_{11}b_{11})c_{12}=\delta (a_{11})b_{11}c_{12}+a_{11}\delta (b_{11})c_{12}.$$
This yields that $\delta (a_{11}b_{11})=\delta
(a_{11})b_{11}+a_{11}\delta (b_{11})$ since $\mathcal {M}$ is a
faithful left $\mathcal {A}$-module.
 \end{proof}

 In the similar manner, one can get
 \begin{lemma}\label{22}
 $\delta $ is a derivation of $\mathcal {T}_{22}$.
 \end{lemma}

Now we can get our main result of this note.

\begin{theorem} \label{theorem} Let $\mathcal {A}$ and $\mathcal {B}$  be two algebras over a $2$-torsion free commutative ring $\mathcal {R}$ with the property:

 (P) Suppose that  $a\in \mathcal {A}$ (resp. $\mathcal {B}$). If $xay+yax=0$ holds for all $x, y\in \mathcal {A}$ (resp. $\mathcal {B}$), then $a=0$.

Let $\mathcal {M}$ be a faithful $(\mathcal {A}, \mathcal
{B})$-bimodule and $\mathcal {T}$ be the triangular algebra
$Tri(\mathcal {A}, \mathcal {M}, \mathcal {B})$. Then every Jordan
derivation $\delta $ on $\mathcal {T}$ into itself is a derivation.
\end{theorem}
\begin{proof}
It follows from Lemmas \ref{1122}-\ref{22}.
\end{proof}

We end this note with the following remark.

\begin{remark} It is easy to see that Theorem \ref{zhangtheorem} is a special case  of Theorem \ref{theorem} when both $\mathcal {A}$ and $\mathcal {B}$ are unital. In other words, our result is a generalization of Theorem \ref{zhangtheorem}.
\end{remark}
 \bibliographystyle{amsplain}

\begin{thebibliography}{10}


\bibitem {ben} D. Benkovi$\check{\textrm{c}}$, Jordan derivations and antiderivations on triangular matrices, \textit{Linear Algebra Appl.} \textbf {397} (2005), 235--244.

\bibitem {bresar} M. Bre$\check{\textrm {s}}$ar, Jordan derivations on semiprime rings, \textit{Proc. Amer. Math. Soc.}, \textbf {104} (1988), 1003--1006.


\bibitem{cheung} W. S. Cheung, Commuting maps on triangular algebras, \textit {J. London Math. Soc.}, \textbf {63} (2001), 117--127.

  \bibitem {ji} P. Ji,   Jordan maps on triangular algebras, \textit {Linear Algebra Appl.}, (to appear).

 \bibitem{ma695} W. S. Martindale III, When are multiplicative mappings additive?\textit {Proc. Amer. Math. Soc.}, \textbf {21} (1969) 695--698.

\bibitem{zhang} J. Zhang, W. Yu, Jordan derivations of triangular algberas, \textit {Linear Algebra Appl.}, \textbf {419} (2006), 251--255.
 \end{thebibliography}

\end{document}